\documentclass{article}
\usepackage{amsfonts}

\begin{document}

\def\Tr{\mbox{Tr \-}}
\def\inj{\mbox{inj \-}}
\def\R{{\mathbb R}}
\def\C{{\mathbb C}}
\def\H{{\mathbb H}}
\def\Ca{{\mathbb Ca}}
\def\Z{{\mathbb Z}}
\def\N{{\mathbb N}}
\def\Q{{\mathbb Q}}
\def\Ad{\mbox{Ad \-}}
\def\k{{\bf k}}
\def\l{{\bf l}}
\def\sp{\mbox{\bf sp}}
\def\Hol{\mbox{Hol}}
\def\Iso{\mbox{Iso}}
\def\Det{\mbox{Det}}
\def\pr{\mbox{pr}}
\def\tr{\mbox{tr}}
\def\sgn{\mbox{sgn}}
\def\const{\mbox{const}}

\title{Globally Hyperbolic Lorentzian Manifolds with Special Holonomy Groups}
\author{Ya.~V.~Baza\u\i{}kin\footnote{The author was supported by the Russian Foundation
for Basic Research (Grant 09--01--00142--a), the President of the
Russian Federation (Grant MK--5482.2008.1), and the Joint Project of
the Siberian and Ural Divisions of the Russian Academy of Sciences
(Grant No.~46).}}
\date{}
\maketitle

%\topmatter
                                          %%% specific SMJ tags
%\translator  A.~P.~Ulyanov\endtranslator
%\iauthor     Baza\u\i{}kin~Ya.~V.\endiauthor
%\UDclass     514.764.214\endUDclass
%\opages      721--736\endopages
                                   %%% standard AMS tags
%\author      Ya.~V.~Baza\u\i{}kin\endauthor
%\title       Globally Hyperbolic Lorentzian Spaces
%             with Special Holonomy Groups
%\endtitle

%\date        March 16, 2009\enddate
%\address     Novosibirsk\endaddress
%\affil       Ya.~V.~Baza\u\i{}kin\\
%             Sobolev Institute of Mathematics,
%             Novosibirsk, Russia
%\endaffil
%\email       bazaikin\@math.nsc.ru
%\endemail
%\thanks{The author was supported by the Russian Foundation
%for Basic Research (Grant 09--01--00142--a), the President of the
%Russian Federation (Grant MK--5482.2008.1), and the Joint Project of
%the Siberian and Ural Divisions of the Russian Academy of Sciences
%(Grant No.~46).}

%\keywords
%holonomy group,
%Lorentzian manifold
%\endkeywords

%\abstract
%We prove that
%each special Lorentzian holonomy group
%(with the exception of those
%including the isotropy groups
%of K\"ahler symmetric spaces)
%can be realized as the holonomy group
%of a~globally hyperbolic Lorentzian manifold.
%\endabstract
%\endtopmatter

\section[]{Introduction}

Well known is the classification of holonomy groups
of the simply connected Riemannian manifolds.
The classical de~Rham decomposition theorem~[1]
immediately reduces the classification problem
to the problem of studying irreducible holonomy groups,
while every irreducible connected Riemannian holonomy group
is either the holonomy group of a~symmetric space
or appears on Berger's list.
Moreover,
each group on Berger's list
is realized as the holonomy group of a~complete Riemannian space.
All relevant references can be found in~[2].

The presence of indecomposable but  irreducible holonomy groups
complicates the situation in the pseudo-Riemannian case. In more
detail, consider a~pseudo-Riemannian manifold $(N,g)$ with the
holonomy group $G=\Hol_p (N)$ for $p \in N$. A~holonomy
representation is called {\it decomposable\/}
whenever there is a~%
$G$-%
invariant decomposition
$
T_pN = W_1 \oplus \dots \oplus W_r
$
with
$r\geq 2$
and
$W_i \neq 0$
for all
$i=1, \dots, r$.
Otherwise,
the representation is called {\it indecomposable}.
A~holonomy representation is called {\it irreducible\/}
whenever there is no nontrivial proper~%
$G$-%
invariant subspace $W \subset T_p N$. The de~Rham decomposition
theorem as generalized to the pseudo-Riemannian case reads [1,\,3]:
Every pseudo-Riemannian manifold with a~decomposable holonomy
representation is locally isometric to the product $ ({\Bbb
R}^{k_1},g_1) \times \dots \times ({\Bbb R}^{k_r},g_r), $ where
$k_i=\dim W_i$ and $\Hol_p(N)=H_1 \times \dots \times H_r$.
If, moreover,~%
$N$
is simply connected and geodesically complete
then
$(N,g)$
is isometric to
$
(N_1, g_1) \times \dots \times
(N_r,g_r),
$
where
$H_i$
is the holonomy group of
$(N_i,g_i)$
for
$i=1,\dots, r$.

Some list of the candidates for irreducible holonomy groups of
pseudo-Rie\-man\-ni\-an manifolds was obtained in [4,\,5], and all
these groups were realized in~[5] as the holonomy groups of
pseudo-Riemannian spaces. Inspecting the list in [4,\,5], we see
that in the Lorentzian case there cannot be any irreducible holonomy
groups but $SO(n+1,1)$. Therefore, the classification problem for
the special holonomy groups of Lorentzian spaces reduces to studying
the indecomposable holonomy representations that are not
irreducible.

The holonomy algebras of the indecomposable Lorentzian manifolds
that are not irreducible
were studied in~[6].
Associated to each algebra
${\bold g} \subset {\bold {so}}(n+1,1)$
of this type
is its {\it orthogonal part\/}
${\bold h} \subset {\bold {so}}(n)$;
if the orthogonal part is given then
there exists exactly four types of algebras~%
${\bold g}$
which can potentially be the holonomy algebras of Lorentzian manifolds.
In the following section
we describe all four types of algebras
as well as the corresponding groups in more detail.

It follows from~[7] that if ${\bold g} \subset {\bold {so}}(n+1,1)$
is the holonomy algebra of an~indecomposable Lorentzian manifold
which is not irreducible then its orthogonal part~%
${\bold h}$ is the holonomy algebra of a~Riemannian manifold. In~[6]
some of these types of algebras were realized, also locally, as the
holonomy algebras of locally defined Lorentzian metrics, and in~[8]
(together with results of [7]) the algebras of all four types were
realized.

However,
the question of the global structure of Lorentzian metrics
with special holonomy
is still not  understood fully.
Moreover,
even the statement of the problem is complicated
by the ambiguity in understanding the meaning of ``completeness''
in Lorentzian geometry.
The problem of constructing globally hyperbolic Lorentzian manifolds
for every special type of holonomy groups
was proposed in~[9].
Briefly speaking,
a~globally hyperbolic Lorentzian space
is a~space possessing a~space-like hypersurface
that is met by each inextensible non-space-like curve
exactly at one point~[10].
This is one of the strongest causality conditions,
and the most useful in mathematical physics.
The special holonomy groups were partially
(namely,
Type~2)
realized in [9]
by globally hyperbolic Lorentzian manifolds.

In this article we continue studying
the problem of constructing globally hyperbolic Lorentzian manifolds
with special holonomy groups
proposed in~[9].
Namely,
the main result of the article is

\vskip0.2cm

{\bf Theorem}. {\it
Take the holonomy group~%
$H$
of a~Riemannian space
whose holonomy representation contains as a~direct factor
no representation of the isotropy of a~K\"ahler symmetric space
of rank greater than~$1$.
Then for every special Lorentzian holonomy group~%
$G$
with the orthogonal part~%
$H$
there exists a~globally hyperbolic Lorentzian manifold
with the holonomy group~%
$G$. }

\vskip0.2cm

Observe that in the cases still to be studied $H=U(n) \times H'$,
the representation of~%
$U(n)$ is not the standard one and can be the holonomy
representation only of the appropriate symmetric space.

The next section is devoted to
the construction of metrics with given holonomy groups,
while in the final section
we study the global causality properties
of the constructed metrics.

\section[]{Construction of Lorentzian Metrics with Special
Holonomy Groups}

Consider a~simply connected time-oriented Lorentzian manifold~%
$N$
of dimension
$n+2$;
thus,
it is a~pseudo-Riemannian space with a~metric~%
$g$ of signature $(n+1,1)$. Take $p \in N$ and the holonomy group
$G=\Hol_p (N) \subset SO(n+1,1)=\Iso (T_p N)$. Owing to the de~Rham
decomposition theorem for pseudo-Riemannian spaces cited above, we
assume henceforth that ~$N$ is indecomposable.

Since the classification of irreducible pseudo-Riemannian manifolds implies that
the irreducible holonomy group of a~Lorentzian manifold
can only be the group
$SO(n+1,1)$,
we assume henceforth that~%
$N$
is not irreducible.
Consequently,
there exists a~proper~$G$-invariant subspace~%
$V$
of
$T_p N$
such that ~$g$
is degenerate on~$V$.
Hence,
the one-dimensional distribution
$L=V \cap V^\perp$
and the
$(n+1)$-%
dimensional distribution
$U=L^\perp \supset L$
arise
which are both $G$-invariant.
It is not difficult to see that
%the metric
$g$
determines on the~%
$n$-%
dimensional space
$\widetilde{U}=U/L$
a~well-defined positive definite inner product,
and~%
$G$
induces the action of some group
$H \subset SO(n)$
on~%
$\widetilde{U}$.
The group~%
$H$
is called the {\it orthogonal part\/} of ~%
$G$.
It follows from~[7] that
if~%
$G$
is the holonomy group of a~Lorentzian manifold
then~%
$H$
is the holonomy group of a~Riemannian manifold;
i.e.,
either it appears on Berger's list,
or it is an~isotropy group of a~symmetric space,
or it is the product of groups of these types.

Consider an~isotropic basis for the tangent space
$T_p N$;
i.e.,
a~basis
in which~
$g$
is defined by the matrix
$$
\left( \begin{array}{ccc} 0 & 0 & 1 \\
0 & E_n & 0 \\
1 & 0 & 0 \end{array} \right).
$$
The following representation of the algebra
${\bold{so}}(n+1, 1)$
results:
$$
{\bold{so}}(n+1, 1) = \left\{  \left(
\begin{array}{ccc} a & X & 0 \\
-Y^T & A & -X^T \\
0 & Y & -a \end{array} \right) \mid  A \in {\bold{so}}(n),\ X,Y \in
{\Bbb R}^n,\ a \in {\Bbb R} \right\}.
$$
Without losing generality, we may assume that~%
$L$
is generated by the first coordinate vector.
Then the Lie algebra of the group
$SO(n+1,1)_L$
preserving~%
$L$
is defined as
$$
{\bold{so}}(n+1,1)_L=
\left\{
\left(
\begin{array}{ccc} a & X & 0 \\
0 & A & -X^T \\
0 & 0 & -a \end{array} \right) \mid a \in {\Bbb R},\ X \in {\Bbb
R}^n,\ A \in {\bold{so}} (n) \right\}.
$$
Some attempt to study the Lie algebras~%
${\bold g}$
corresponding to the possible holonomy groups
$G \subset SO (n+1, 1)$
was made in~[6].
Take the Lie algebra~%
${\bold h}$
of the orthogonal part~%
$H$
of a~group
$G \subset SO(n+1,1)_L$.
It is proved in~[6] that~%the algebra
${\bold g}$
can only be of one of the following four types:
$$
{\bold g}^{1,{\bold h}} = \left\{ \left(
\begin{array}{ccc} a & X & 0 \\
0 & A & -X^T \\
0 & 0 & -a
\end{array}
\right) \mid  a \in {\Bbb R},\ X \in {\Bbb R}^n,\ A \in {\bold h} \subset {\bold
{so}}(n) \right\};
$$
$$
{\bold g}^{2,{\bold h}} = \left\{ \left(
\begin{array}{ccc} 0 & X & 0 \\
0 & A & -X^T \\
0 & 0 & 0
\end{array}
\right) \mid  X \in {\Bbb R}^n,\ A \in {\bold h} \subset {\bold{so}}(n)
\right\};
$$
$$
{\bold g}^{3,{\bold h},\phi } = \left\{ \left(
\begin{array}{ccc}
\phi(A) & X & 0 \\
0 & A & -X^T \\
0 & 0 & -\phi(A)
\end{array}
\right) \mid  X \in {\Bbb R}^n,\ A \in {\bold h} \subset {\bold{so}}(n)
\right\},
$$
where the center
$Z({\bold h})$
of~%
${\bold h}$
is nontrivial
and
$\phi: {\bold h} \rightarrow {\Bbb R}$
is a~nonzero linear mapping with
$\phi|_{{\bold h}'}=0$
(we denote by
${\bold h}'$
the commutant of~
${\bold h}$);
$$
{\bold g}^{4,{\bold h}, m, \psi }
= \left\{ \left(
\begin{array}{cccc} 0 & X & \psi(A) & 0 \\
0 & A & 0 & -X^T \\
0 & 0 & 0 & -\psi(A)^T \\
0 & 0 & 0 & 0
\end{array}
\right) \mid  X \in {\Bbb R}^m,\ A \in {\bold h} \subset {\bold{so}}(m)
\right\}
$$
for
$0<m<n$,
where
$\dim Z({\bold h}) \geq n-m$
and
$\psi: {\bold h} \rightarrow {\Bbb R}^{n-m}$
is a~surjective linear mapping with
$\psi|_{{\bold h}'}=0$.

Take the center
$T^r \subset H$
of~%
$H$. Denote by $\Det : H \rightarrow T^r$ the uniquely defined
homomorphism such that $\Det^{-1}(1) \subset H$
is the semisimple part of~%
$H$.
The algebras above are
the tangent Lie algebras to the following subgroups of
$SO(n+1, 1)$:
$$
G^{1,H} = \left\{ \left(
\begin{array}{ccc} e^a & X & -\frac{1}{2} e^{-a} X X^T \\
0 & A & -e^{-a} A X^T \\
0 & 0 & e^{-a}
\end{array}
\right) \mid  a \in {\Bbb R},\ X \in {\Bbb R}^n,\ A \in H \subset SO(n) \right\};
$$
$$
G^{2,H} = \left\{ \left(
\begin{array}{ccc} 1 & X & -\frac{1}{2} X X^T \\
0 & A & - A X^T \\
0 & 0 & 1
\end{array}
\right) \mid  X \in {\Bbb R}^n,\ A \in H \right\};
$$
$$
%\gathered
G^{3,H,\phi} = \left\{ \left(
\begin{array}{ccc} e^{\phi(a_1,\dots, a_r)} & X & -\frac{1}{2} e^{-a} X X^T \\
0 & A & -e^{-a} A X^T \\
0 & 0 & e^{-\phi(a_1,\dots, a_r)}
\end{array}
\right) \mid \right.
$$
$$
\left.
X \in {\Bbb R}^n,\ A \in H, %\right.
%\\
%\left.
%\vphantom{
%\matrix e^{\phi(a_1,\dots, a_r)} & X & -\frac{1}{2} e^{-a} X X^T \\
%0 & A & -e^{-a} A X^T \\
%0 & 0 & e^{-\phi(a_1,\dots, a_r)}
%\endmatrix
%}
\ \Det (A) = (e^{i a_1}, \dots, e^{i a_r}) \in T^r \right\};
%\endgathered
$$
$$
%\gathered
G^{4,H,m,\psi} = \left\{ \left(
\begin{array}{cccc} 1 & X & \psi (a_1, \dots, a_r) & -\frac{1}{2} (X X^T +Y Y^T) \\
0 & A & 0 & - A X^T \\
0 & 0 & E_m & - \psi (a_1, \dots, a_r)^T \\
0 & 0 & 0 & 1
\end{array}
\right) \right.
\\
\left.
\vphantom{
\begin{array}{cccc} 1 & X & \psi (a_1, \dots, a_r) & -\frac{1}{2} (X X^T +Y Y^T) \\
0 & A & 0 & - A X^T \\
0 & 0 & E_m & - \psi (a_1, \dots, a_r)^T \\
0 & 0 & 0 & 1
\end{array}
} \mid X \right.
$$
$$
\left. \in {\Bbb R}^{n-m},\ A \in H \subset SO(n-m),\ \Det (A)=(e^{i
a_1}, \dots, e^{i a_r}) \right\}.
%\endgathered
$$
Observe~[6] that
$G^{3,H,\phi}$
and
$G^{4,H,m,\psi}$
cannot be closed subgroups in
$SO(n+1, 1)$.

Take some Riemannian manifold~%
$M$
of dimension~%
$n$
with metric~%
$g$.
Equip
$N=M \times
{\Bbb R}^2$
with the Lorentzian metric
$$
\tilde{g}= 2 d \eta (d \xi  + \varepsilon f d \eta + 2 \varepsilon A
) + g
                                \eqno{(1)}
$$
for the  coordinates
$\xi, \eta$
on the plane
${\Bbb R}^2$,
some function~%
$f$
on~%
$N$,
some~$1$-%
form~%
$A$
on~%
$M$,
and some real parameter
$\varepsilon>0$.

\vskip0.2cm

{\bf Theorem 1}. {\it
Take the holonomy group~%
$H$ of a~Riemannian space whose holonomy representation contains as
a~direct factor no isotropy group of a~K\"ahler symmetric space of
rank greater than~$1$. For each of the groups $G^{1, H}$, $G^{2,H}$,
$G^{3,H,\phi}$, and $G^{4,H,m,\psi}$ there exists a~Lorentzian
manifold with metric of the form~{\rm (1)} realizing this holonomy
group. }

\vskip0.2cm

We devote the remainder of this section
to proving the theorem.
Take an~orthonormal coframe
$e^1, \dots, e^n$
of the metric~%
$g$ which is only locally defined in general. Extend it to
an~isotropic coframe $\tilde{e}^0, \tilde{e}^1, \dots, \tilde{e}^n,
\tilde{e}^{n+1}$
of~%
$\tilde{g}$
as follows:
$$
\tilde{e}^0= d\xi + \varepsilon f d\eta +2 \varepsilon A,
\quad
\tilde{e}^i=e^i,
\quad
\tilde{e}^{n+1}= d\eta.
$$
Henceforth we agree that
the Greek indices will take values from~%
$0$
to
$n+1$
and the Latin indices,
from~%
$1$
to~%
$n$.
Take the dual frame
$\tilde{e}_\alpha$
to
$\tilde{e}^\alpha$.
Recall that the connection form
$\widetilde{\omega}$
and the curvature form
$\widetilde{\Omega}$
are found from the relations
$d \tilde{e}^\alpha =
-\widetilde{\omega}^\alpha_\beta \wedge \tilde{e}^\beta$
and
$\widetilde{\Omega}^\alpha_\beta=d \widetilde{\omega}^\alpha_\beta+
\widetilde{\omega}^\alpha_\gamma \wedge \widetilde{\omega}^\gamma_\beta$,
where the matrices
$\bigl(\widetilde{\omega}^\alpha_\beta\bigr)_{\alpha,\beta}$
and
$\bigl(\widetilde{\Omega}^\alpha_\beta\bigr)_{\alpha,\beta}$
lie in the algebra
${\bold {so}}(n+1,1)$.

Direct calculations
which we omit
lead us to the following statement:

\vskip0.2cm

{\bf Lemma 1}. {\it In the isotropic coframe $\tilde{e}^0,
\tilde{e}^1, \dots, \tilde{e}^n, \tilde{e}^{n+1}$ the torsion and
curvature forms of the metric~ $\tilde{g}$ are as follows:
$$
%\gathered
\widetilde{\omega}_0^0  =  -\widetilde{\omega}_{n+1}^{n+1}  =
\varepsilon f_0 \tilde{e}^{n+1},
\quad
\widetilde{\omega}_{n+1}^0  =  -\widetilde{\omega}_0^{n+1}  =  0,
\quad
\widetilde{\omega}_i^{n+1}  =  -\widetilde{\omega}_0^i  =  0,
$$
$$
\widetilde{\omega}_i^0  =  -\widetilde{\omega}_{n+1}^i  =
\varepsilon f_i \tilde{e}^{n+1}+ \varepsilon F_{ij} \tilde{e}^j,
\quad \widetilde{\omega}_j^i  =  -\widetilde{\omega}_i^j  =
\omega_j^i- \varepsilon F_{ij} \tilde{e}^{n+1},
%\endgathered
%\tag2
$$
where
$f_0=\langle d f, e^{n+1} \rangle = \partial f/\partial \xi$,
$f_i=\langle d f, e^i \rangle$,
and
$d A=F=\frac{1}{2} F_{ij} e^i\wedge e^j$,
$$
%\gather
\widetilde{\Omega}_0^0  =  -\widetilde{\Omega}_{n+1}^{n+1}  =
\varepsilon f_{00} \tilde{e}^0 \wedge \tilde{e}^{n+1}+ \varepsilon f_{0i} \tilde{e}^i \wedge \tilde{e}^{n+1},
\\
\widetilde{\Omega}_0^{n+1}  =  -\widetilde{\Omega}_{n+1}^0  =  0,
$$
$$
\widetilde{\Omega}_0^i  =  -\widetilde{\Omega}_i^{n+1}  =  0,
%\tag3
\widetilde{\Omega}_i^0 =  -\widetilde{\Omega}_{n+1}^i  =
$$
$$
\varepsilon \widetilde{\nabla}_0 f_i \tilde{e}^0 \wedge
\tilde{e}^{n+1} + (\varepsilon \widetilde{\nabla}_k f_i +
\varepsilon^2 F_{ij} F_{jk}) \tilde{e}^k \wedge \tilde{e}^{n+1} +
\varepsilon \widetilde{\nabla}_k F_{ij} \tilde{e}^k \wedge
\tilde{e}^{j},
$$
$$
\widetilde{\Omega}_j^i  =  -\widetilde{\Omega}_i^j  =  \Omega_j^i -
\varepsilon \nabla_k F_{ij} \tilde{e}^k \wedge \tilde{e}^{n+1},
%\endgather
$$
where $f_{00}=\partial^2 f /\partial \xi^2$ and $f_{0i}=\langle f_0,
e^i\rangle$. }

\vskip0.2cm

In accordance with the Ambrose--Singer theorem~ [11]
the holonomy algebra
${\bold {hol}}={\bold {hol}}_p (N)$
of~%manifold
$N$
is generated by the elements of the form
$
(P_\gamma \widetilde{\Omega} ) (v,w) \in {\bold{so}}(n+1, 1),
$
where
$v,w \in T_p N$
and~%
$\gamma$
is a~path in~%
$N$
ending at ~%
$p$,
while
$P_\gamma$
stands for the parallel transport along~
$\gamma$.
Since parallel transport preserves the isotropy
of coframes,
at every point we can identify
the group of isometries of the tangent space
with the matrix group
$SO(n+1,1)$
preserving the inner product~(1).

Take a~smooth curve
$\gamma(t)$
with
$\gamma (0)=p$
and transport the coframe
$\tilde{e}^\alpha$
along~%
$\gamma$:
$$
(P_\gamma)_t (\tilde{e}^\alpha)= (W_t)^\alpha_\beta \tilde{e}^\beta,
\quad
W_t \in SO (n+1, 1).
$$
The Ambrose--Singer theorem implies that
the algebra
${\bold{hol} }$
is generated by the elements
$Ad (W_t) \widetilde{\Omega} (v,w)$
for all
$p\in N$,
$v, w \in T_p N$,
and all sufficiently small~%
$t$.
Passing to the limit as
$t \rightarrow 0$,
we find that the algebra
${\bold{hol} }$
is generated by the matrices
$\widetilde{\Omega}(v,w)$
and
$[ \widetilde{\omega}(u), \widetilde{\Omega}(v,w) ]$
for all
$p \in N$
and
$u,v,w \in T_p N$.
Therefore,
we have proved the following statement:

\vskip0.2cm

{\bf Lemma 2}. {\it The holonomy algebra ${\bold{hol}} \subset
{\bold {so}}(n+1,1)$
of ~%
$N$
can be found as
$$
{\bold{hol}}= {\cal L}_{\Bbb R} \{
\widetilde{\Omega}(\tilde{v}_\alpha \wedge \tilde{v}_\beta),
[\widetilde{\omega}(\tilde{v}_\gamma),
\widetilde{\Omega}(\tilde{v}_\alpha \wedge \tilde{v}_\beta)] \mid  p
\in N, \alpha, \beta, \gamma = 0, 1, \dots, n, n+1 \}.
$$
}

\vskip0.2cm

Take an~algebra
${\bold g} \subset {\bold{so}}(1,n+1)$
of one the four types defined above,
with the orthogonal part~%
${\bold h}$.
The argument above implies that~%
${\bold h}$
is the holonomy algebra of a~Riemannian manifold.
Using the classification of Riemannian holonomy groups,
we obtain the orthogonal decomposition
$$
\widetilde{U}_p= {\Bbb R}^n= {\Bbb R}^{n-m} \oplus {\Bbb R}^{n_0} \bigoplus\limits_{i=1}^r {\Bbb R}^{n_i},
\quad
m=n_0+\sum\limits_{i=1}^r n_i,
$$
as well as the corresponding decomposition
$$
{\bold h} = {\bold 0} \oplus {\bold h}_0 \bigoplus\limits_{i=1}^r {\bold h}_i \subset
{\bold{so}}(n)
$$
of algebras,
where
${\bold 0} \subset {\bold{so}}(n-m)$
is the trivial term,
%the algebra
${\bold h}_0$
is a~
(possibly reducible) subalgebra
of
${\bold {so}}(n_0)$
with the trivial center,
while every algebra
${\bold h}_i$
for
$i=1,
\dots, r$
is isomorphic to
${\bold u}(m_i)$,
where
$2 m_i=n_i$,
with the standard action on
${\Bbb R}^{n_i}=\Bbb C^{m_i}$,
meaning that
${\bold h}_i ({\Bbb R}^{n_j})=0$
for
$i,j=0, \dots, r$
with
$i\neq j$
and
${\bold h}$
annihilates the term
${\Bbb R}^{n-m}$.

The following notation will be convenient:
$$
{\cal A}= \left\{ \left(
\begin{array}{ccc} a & 0 & 0 \\
0 & 0 & 0 \\
0 & 0 & -a
\end{array}
\right) \mid  a \in {\Bbb R} \right\}, \quad {\cal K}= \left\{
\left(
\begin{array}{ccc} 0 & 0 & 0 \\
0 & A & 0 \\
0 & 0 & 0
\end{array}
\right) \mid  A \in {\bold{so}}(n) \right\},
$$
$$
{\cal N}= \left\{ \left(
\begin{array}{ccc} 0 & X & 0 \\
0 & 0 & -X^T \\
0 & 0 & 0
\end{array}
\right) \mid  X \in {\Bbb R}^n \right\}.
$$
Therefore, ${\bold{so}}(n+1,1)_L={\cal A} \oplus {\cal K} \oplus
{\cal N}$
for an~abelian ideal~%
${\cal N}$ in ${\bold{so}}(n+1,1)_L$
and a~subalgebra~%
${\cal K}$ isomorphic to ${\bold{so}}(n)$
and commuting with~%
${\cal A}$. Denote the elements of ${\bold{so}}(n+1,1)_L$ by $(a, A,
X)$ in accordance with the decomposition. Finally, denote the
projections of ${\bold{so}}(n+1,1)_L$ onto the subalgebras of our
decomposition by $\pr_{\cal A}$, $\pr_{\cal K}$, and $\pr_{\cal N}$.

Consider these types of algebra by cases.

\medskip
{\bf Type 1.}
In this case take an~
$n$-dimensional compact Riemannian manifold~%
$M$
with the holonomy algebra~%
${\bold h}$,
put
$N=M
\times {\Bbb R}^2$,
and consider on~%
$N$
the metric~(1) with
$A=0$.
It is clear from~(3) that
if the function~%
$f$ is of a~sufficiently general form then $ {\cal L}_{\Bbb R} \{
\widetilde{\Omega}(\tilde{v}_\alpha \wedge \tilde{v}_\beta) \mid
\alpha, \beta \} = {\bold g}^{1, {\bold h}}. $
Furthermore,~%
(2) shows that
$\widetilde{\omega}
(\tilde{v}_\gamma) \in {\bold g}^{1, {\bold h}}$,
which yields the claim.

\medskip
{\bf Type 2}.
This case is completely analogous to the previous,
as we only have to consider the function~%
$f$
of a~sufficiently general form
and independent of ~%
$\xi$
and~%
$\eta$.

In order to continue the proof
we will need the Calabi example of a~Riemannian space
with the holonomy group
$SU(n)$.
Our exposition of this construction follows~[12].
Take the complex projective space
$\Bbb C P^{n-1}$
with the Fubini--Study metric
$ds^2$.
Consider the metric
$$
d \hat{s}^2 = \frac{d \rho^2}{1-\frac{1}{\rho^{2n}}} + \rho^2 \left(
1-\frac{1}{\rho^{2n}} \right) (d\tau - 2 A)^2 + \rho^2 ds^2,
                                \eqno{(4)}
$$
where~%
$A$
is a~%
$1$-%
form on
$\Bbb C P^{n-1}$
such that
$d A=\Phi$
is a~K\"ahler form,
and
$\rho \geq 1$,~%
$\tau$
are new variables,
while~%
$\tau$
is periodic.
The form~%
$A$
is only defined  locally on
${\Bbb C}P^{n-1}$;
however,
since
${\Bbb C}P^{n-1}$
is a~Hogde manifold,
we can choose the period
$\Delta \tau$
of~%
$\tau$
so that
$\int 2 \Phi$
over every closed~%
$2$-%
chain is an~integer multiple of
$\Delta \tau$.
Then
$d\tau-2 A$
is independent of the choice of a~coordinate neighborhood,
and~(4) is a~well-defined smooth metric
on the total space of the complex line bundle over
${\Bbb C}P^{n-1}$
which is the~%
$n$th power of the Hopf bundle.
%${\Bbb C}P^{n-1}$.
Moreover,
the resulting complete Riemannian manifold
$C_n$
(the Calabi space)
possesses the holonomy group
$SU(n)$.
The space
$C_n$,
constructed by Calabi~[13],
is a~generalization of the Eguchi--Hanson space~[14]
arising for
$n=2$.
Moreover,
the K\"ahler form on
$C_n$
is
$
\widehat{\Phi}=- \rho^2 \Phi + \rho d\rho \wedge (d\tau - 2 A).
$

Below we will need the~%
$1$-%
form
$
B=\frac{1}{2} \rho^2 (d\tau - 2 A).
$
By above,~%
$B$
is globally defined on~
$C_n$
and
$d B = \widehat{\Phi}$.

\medskip
{\bf Type 3}.
Take a~compact Riemannian manifold
$M^{n_0}$
with the holonomy algebra
${\bold h}_0$,
and take the direct product
$M^n=M^{n_0} \times C_{m_1}
\times \dots \times C_{m_r}$
of Riemannian spaces,
with
$n= n_0+ 2\sum\nolimits_{i=1}^r m_i$,
where the Calabi space
$C_{m_i}$
is defined above.
Denote by
$\rho_i$
and
$\tau_i$
the corresponding coordinates,
while by
$B_i$,
for
$i=1, \dots, r$,
the~%
$1$-%
forms on
$C_{m_i}$.
It is not difficult to observe that
$\rho_i^2$
and
$B_i$
are globally defined and smooth on the whole of~%
$N$.
The center of the algebra
$Z({\bold h})=\bigoplus\nolimits_{i=1}^r {\bold h}_i$
is isomorphic to
$\bigoplus\nolimits_{i=1}^r {\bold u}(1)$,
and we can choose as a~basis for
$Z({\bold h})$
the forms
$\widehat{\Phi}_i$
for
$i=1,\dots, r$
(it is natural to interpret~%
$2$-%
forms as the elements of~%
${\cal K}$). Then a~collection of real constants $\phi_i$, for $i=1,
\dots, r$, not vanishing simultaneously, determines a~nonzero linear
mapping $\phi: Z({\bold h}) \rightarrow {\Bbb R}$. Put $ A=
\sum\nolimits_{i=1}^r g_i B_i, $ where each function $g_i$ depends
only on the variable $\rho_i$. Take
$$
f=h - \xi \sum\limits_{i=1}^r \phi_i \left(\frac{1}{2m_i}\rho_i g_i' + g_i
\right),
$$
where the function~%
$h$
is independent of~%
$\xi$
and~
$\eta$.
Put
$N=M \times {\Bbb R}^2$
and equip it
with the metric~(1)
defined by ~%
$f$
and ~%
$A$.

Put
$$
W_1={\cal L}_{\Bbb R} \{ \widetilde{\omega}(\tilde{e}_k) \mid  k=1,
\dots, n \}, \quad W_2={\Bbb R} \widetilde{\omega}(\tilde{e}_{n+1}).
$$
It is clear that $ \pr_{\cal K} (W_1)={\cal L}_{\Bbb R} \{
\omega(e_k) \mid k=1,\dots, n \} $ and $ \pr_{\cal A}(W_1)=0$.
Furthermore,
$$
%\gathered
F=dA=\sum\limits_{i=1}^r
(dg_i \wedge B_i + g_i d B_i)
%\\
= \sum\limits_{i=1}^r \bigl(g_i' d\rho_i \wedge B_i
+ g_i \widehat{\Phi}_i\bigr)
=
$$
$$
\sum\limits_{i=1}^r \left( \frac{1}{2}\rho_i^2 g_i' d \rho_i \wedge
(d \tau_i - A_i)+ g_i \widehat{\Phi}_i\right).
%\endgathered
$$
Therefore, $F \in\bigoplus\nolimits_{i=1}^r{\bold h}_i$. By (2) we
have $\pr_{\cal K}(\widetilde{\omega}(\tilde{e}_{n+1}))=-
\varepsilon F$. Consequently,
$$
\pr_{{\bold{so}}(n_0)}(\widetilde{\omega}(\tilde{e}_{n+1}))=0, \quad
\pr_{{\bold h}_i}(\widetilde{\omega}(\tilde{e}_{n+1}))=-\varepsilon
\left(\frac{1}{2}\rho_i^2 g_i' d \rho_i \wedge (d \tau_i - A_i)+ g_i
\widehat{\Phi}_i\right).
$$
Every algebra
${\bold h}_i$
is isomorphic to
${\bold u}(m_i)$;
therefore,
the projection
$
\tr : {\bold h}_i \rightarrow Z({\bold h}_i)={\bold u}(1)
$
is defined
which is invariant under the holonomy group.
Moreover,
the form
$\widehat{\Phi}_i$
is a~generator for the center
$Z({\bold h}_i)$,
so that
$$
\pr_{Z({\bold h}_i)}(W_2)=-\varepsilon \tr \left(
\frac{1}{2}\rho_i^2 g_i' d \rho_i \wedge (d \tau_i - A_i)+ g_i
\widehat{\Phi}_i \right)= -\varepsilon \left( \frac{1}{2m_i}\rho_i
g_i' + g_i \right) \widehat{\Phi}_i.
$$
Now,
$$
\pr_{\cal A} (\widetilde{\omega}(\tilde{e}_{n+1})) = -\varepsilon
\sum\limits_{i=1}^r \phi_i \left(\frac{1}{2m_i}\rho_i g_i' + g_i
\right)= \phi ( \pr_{Z({\bold
h})}(\widetilde{\omega}(\tilde{e}_{n+1})) ).
$$
Therefore, ${\cal L}_{\Bbb R} \{\widetilde{\omega}(\tilde{e}_\alpha)
\mid \alpha \}=W_1\oplus W_2 \subset {\bold g}^{3,{\bold h}, \phi}$.

Now, we put
$$
V_1 = {\cal L}_{\Bbb R} \{ \widetilde{\Omega} (\tilde{e}_i \wedge
\tilde{e}_j) \mid  i,j =1, \dots, n \} \subset {\bold{so}} (n+1,1).
$$
Inspecting~(3) and applying Lemma~2,
we see that
$$
\pr_{\cal
 K}(V_1) = {\cal L}_{\Bbb R}
\{ \Omega(e_i \wedge e_j)\mid  i,j=1, \dots, n \} \subset
{\bold h}_0 \bigoplus\limits_{i=1}^r {\bold h}_i',
$$
where ${\bold h}_i'={\bold {su}}(m_i)$. Moreover, $\pr_{\cal A}
(V_1)=0$. Furthermore, put
$$
V_2 = {\cal L}_{\Bbb R} \{ \widetilde{\Omega} (\tilde{e}_k \wedge
\tilde{e}_{n+1}) \mid  k=1, \dots, n \}.
$$
Since
$$
F=\sum\limits_{i=1}^r \left( \frac{1}{2}\rho_i^2 g_i' d \rho_i \wedge
(d \tau_i - A_i)+ g_i \widehat{\Phi}_i\right) \in \bigoplus\limits_{i=1}^r{\bold h}_i,
$$
we have $\nabla_k F \in \bigoplus\nolimits_{i=1}^r{\bold h}_i$.
By~(3), $\pr_{\cal K}(\widetilde{\Omega}(\tilde{e}_k\wedge
\tilde{e}_{n+1}))=-\varepsilon \nabla_k F$. Consequently,
$$
%\gathered
\pr_{{\bold{so}}(n_0)}(\widetilde{\Omega}(\tilde{e}_k\wedge
\tilde{e}_{n+1}))=0,
\\
\pr_{{\bold h}_i}(\widetilde{\Omega}(\tilde{e}_k \wedge
\tilde{e}_{n+1}))=
$$
$$
-\varepsilon \nabla_k \left( \frac{1}{2}\rho_i^2 g_i' d \rho_i
\wedge (d \tau_i - A_i)+ g_i \widehat{\Phi}_i\right)
%\endgathered
$$
for
$i=1,\dots, r$
and
$k=1, \dots, n$.
Then,
since the form
$\widehat{\Phi}_i$
is parallel,
$$
%\allowdisplaybreaks\gathered
\pr_{Z({\bold h}_i)}(\widetilde{\Omega}(\tilde{e}_k\wedge
\tilde{e}_{n+1}))=-\varepsilon \tr \left( \nabla_k \left(
\frac{1}{2}\rho_i^2 g_i' d \rho_i \wedge (d \psi_i - A_i)+ g_i
\widehat{\Phi}_i \right) \right) =
$$
$$
-\varepsilon \nabla_k \tr \left( \frac{1}{2}\rho_i^2 g_i' d \rho_i
\wedge (d \psi_i - A_i)+ g_i \widehat{\Phi}_i \right) = -\varepsilon
\nabla_k \left( \left( \frac{1}{2m_i}\rho_i g_i' + g_i \right)
\widehat{\Phi}_i \right) =
$$
$$
 -\varepsilon
\left(1-\frac{1}{\rho_i^{2m_i}}\right)\left(\frac{2 m_i +1}{2 m_i}
g_i' + \frac{\rho_i}{2 m_i} g_i''\right) \widehat{\Phi}_i
%\endgathered
$$
if the index~%
$k$ corresponds to the variable $\rho_i$ and $ \pr_{Z({\bold
h}_i)}(\widetilde{\Omega}(\tilde{e}_k\wedge \tilde{e}_{n+1}))=0 $
otherwise. Furthermore,
$$
\pr_{\cal A} (\widetilde{\Omega}(\tilde{e}_k \wedge
\tilde{e}_{n+1})) = -\varepsilon \phi_i \left(
1-\frac{1}{\rho_i^{m_i}} \right) \left(\frac{1}{2m_i}\rho_i g_i' +
g_i \right)^\prime = \phi ( \pr_{Z({\bold
h})}(\widetilde{\Omega}(\tilde{v}_i \wedge \tilde{v}_{n+1})))
$$
in the case that~% the index~%
$k$
corresponds to %the variable
$\rho_i$.
Finally,
$V_3={\Bbb R}\widetilde{\Omega}(\tilde{e}_0 \wedge
\tilde{e}_{n+1})=0$.
The last relations yield
$$
{\cal L}_{\Bbb R} \{ \widetilde{\Omega} (\tilde{v}_\alpha,
\tilde{v}_\beta) \mid \alpha, \beta \} = V_1+V_2+V_3 \subset {\bold
g}^{3,{\bold h},\phi}.
$$
Moreover,
$$
\pr_{\cal K} ([W_1, V_1] \oplus V_1)={\bold h_0}
\bigoplus\limits_{i=1}^r {\bold h}_i'
$$
by Lemma~2 applied to the manifold~%
$M$.
It remains to observe that
if we choose a~sufficiently general function~%
$h$ then $\pr_{\cal N}(V_2)={\cal N}$, which yields
$$
{\cal L}_{\Bbb R} \{
[\widetilde{\omega}(\tilde{e}_\gamma),\widetilde{\Omega}
(\tilde{v}_\alpha, \tilde{v}_\beta) ], \widetilde{\Omega}
(\tilde{v}_\alpha, \tilde{v}_\beta) \mid \alpha, \beta\, \gamma \} =
{\bold g}^{3,{\bold h},\phi}.
$$
Lemma~2 implies that
${\bold{hol}}(N) = {\bold g}^{3,{\bold h},\phi}$.
Since the orthogonal coframe
$e^1, \dots, e^n$
is chosen in an~arbitrary neighborhood of~%
$M$, this implies that $\Hol (N)=G^{3,H,\phi}$.

\medskip
{\bf Type 4}.
The proof is generally similar to the previous case.
Consider a~compact Riemannian manifold
$M^{n_0}$
with the holonomy algebra
${\bold h}_0$.
Take the direct product
$$
M^n={\Bbb R}^{n-m} \times M^{n_0} \times C_{m_1} \times \dots \times C_{m_r},
\quad
m= n_0+ 2\sum\limits_{i=1}^r m_i,
$$
of Riemannian spaces,
where
${\Bbb R}^{n-m}$
is the
$(n-m)$-%
dimensional flat Euclidean space with the variables
$t_1,\dots, t_{n-m}$,
while the Calabi space
$C_{m_i}$
is defined above.
Denote by
$\rho_i$,
$\tau_i$,
and
$B_i$,
for
$i=1, \dots, r$,
the corresponding coordinates and~%
$1$-%
forms on
$C_{m_i}$.
As above,
$\rho_i^2$
and
$B_i$
are globally defined and smooth on the whole of~%
$N$.
Define the linear mapping
$\psi: Z({\bold h})=\bigoplus\nolimits_{i=1}^r Z({\bold h}_i)
\rightarrow {\Bbb R}^{n-m}$
in the basis
$\widehat{\Phi}_i$,
for
$i=1, \dots, r$,
by the matrix
$(\psi_i^j)_{i=1}^r;_{j=1}^{n-m}$
of maximal rank.

Put
$
A= \sum\nolimits_{i=1}^r g_i B_i,
$
where each function
$g_i$
depends only on the variable
$\rho_i$,
and
$$
f=h- \sum\limits_{k=1}^{n-m} \sum\limits_{i=1}^r \psi_i^k t_k
\left(\frac{1}{2m_i}\rho_i g_i' + g_i \right)
$$
for some function~%
$h$
independent of~%
$\xi$
and~%
$\eta$.
Put
$N=M
\times {\Bbb R}^2$.
Consider on~%
$N$
the metric~(1) defined by the function~%
$f$
and the~%
$1$-%
form~%
$A$.

\goodbreak

Put $W_1={\cal L}_{\Bbb R} \{ \widetilde{\omega}(\tilde{e}_k) \mid
k=1, \dots, n \}$ and $W_2={\Bbb R}
\widetilde{\omega}(\tilde{e}_{n+1})$. As above,
$$
%\gather
\pr_{\cal K} (W_1)={\cal L}_{\Bbb R} \{ \omega(e_k) \mid k=1, \dots,
n \}, \quad \pr_{\cal K}(W_2) \subset {\bold h},
$$
$$
\pr_{Z({\bold h}_i)}(\widetilde{\omega}(\tilde{e}_{n+1}))=
-\varepsilon \left( \frac{1}{2m_i}\rho_i g_i' + g_i \right)
\widehat{\Phi}_i.
%\endgather
$$
Denote by ${\cal N}_1$ and ${\cal N}_2$
the subspaces of~%
${\cal N}$ corresponding to ${\Bbb R}^{n-m}$ and $M^{n_0}\times
C_{m_1} \times \dots \times C_{m_r}$ respectively. It is clear that
${\cal N}={\cal N}_1 \oplus {\cal N}_2$. Furthermore, $\pr_{{\cal
N}_1} (W_1)=0$ and
$$%\align
\pr_{{\cal N}_1} (\widetilde{\omega}(\tilde{e}_{n+1}))  = -
\varepsilon \sum\limits_{k=1}^{n-m} \frac{\partial f}{\partial
t_k}\frac{\partial}{\partial t_k}
\\
= -\varepsilon \sum\limits_{k=1}^{n-m} \sum\limits_{i=1}^r \psi_i^k
\left(\frac{1}{2m_i}\rho_i g_i' + g_i \right)
\frac{\partial}{\partial t_k}=
$$
$$
\psi (\pr_{Z({\bold h})} (\widetilde{\omega}(\tilde{e}_{n+1}))),
%\endalign
$$
where $\frac{\partial}{\partial t_k}$, for $k=1, \dots, n-m$, is
a~basis for ${\Bbb R}^{n-m}$. Therefore, ${\cal L}_{\Bbb R} \{
\widetilde{\omega}(\tilde{e}_\alpha) \mid\alpha \} \subset {\bold
g}^{4,{\bold h}, m, \psi}$. As above, put
$$
V_1 = {\cal L}_{\Bbb R} \{ \widetilde{\Omega} (\tilde{e}_i \wedge
\tilde{e}_j) \mid  i,j =1, \dots, n \}.
$$
From~(3) we see that
$$
\pr_{{\bold{so}}(n)}(V_1) \subset {\bold h}_0
\bigoplus\limits_{i=1}^r {\bold h}_i'
$$
and $\pr_{{\cal N}_1}(V_1)=0$.
Furthermore,~%
%the function
$f$
is independent of~%
$\xi$;
therefore,
the space
$V_2 =
{\Bbb R}\widetilde{\Omega} (\tilde{e}_0 \wedge \tilde{e}_{n+1})$
is trivial.
Finally,
consider the space
$$
V_3 = {\cal L}_{\Bbb R} \{ \widetilde{\Omega} (\tilde{e}_i \wedge
\tilde{e}_{n+1}) \mid  i =1, \dots, n \}.
$$
We have $\pr_{\cal K}(\widetilde{\Omega}(\tilde{e}_k\wedge
\tilde{e}_{n+1}))=-\varepsilon \nabla_k F$; consequently,
$\pr_{{\bold{so}}(n-m+n_0)}(\widetilde{\Omega}(\tilde{e}_k\wedge
\tilde{e}_{n+1}))=0$ and
$$
\pr_{{\bold h}_i}(\widetilde{\Omega}(\tilde{e}_k\wedge
\tilde{e}_{n+1}))=-\varepsilon \nabla_k \left( \frac{1}{2}\rho_i^2
g_i' d \rho_i \wedge (d \psi_i - A_i)+ g_i \widehat{\Phi}_i\right)
$$
for
$i=1,\dots, r$
and
$k=1, \dots, n$.
Then
$$
\pr_{Z({\bold h}_i)}(\widetilde{\Omega}(\tilde{e}_k\wedge
\tilde{e}_{n+1}))=-\varepsilon
\biggl(1-\frac{1}{\rho_i^{2m_i}}\biggr) \left(\frac{2 m_i +1}{2 m_i}
g_i' + \frac{\rho_i}{2 m_i} g_i''\right) \widehat{\Phi}_i
$$
if the index~%
$k$ corresponds to the variable $\rho_i$, and $ \pr_{Z({\bold
h}_i)}(\widetilde{\Omega}(\tilde{e}_k\wedge \tilde{e}_{n+1}))=0 $
otherwise. If~$k$ corresponds to $\rho_i$ then~(3) implies that
$$
%\gathered
\pr_{{\cal
N}_1}(\widetilde{\Omega}(\tilde{e}_k\wedge\tilde{e}_{n+1}))=-
\varepsilon \sum\limits_{j=1}^{n-m} \widetilde{\nabla}_k f_j
\frac{\partial}{\partial t_j}= -\varepsilon \sum\limits_{j=1}^{n-m}
\nabla_k \sum\limits_{i=1}^r \psi_i^j \biggl(\frac{1}{2m_i}\rho_i
g_i' + g_i \biggr) \frac{\partial}{\partial t_j}
$$
$$
= -\varepsilon \sum\limits_{j=1}^{n-m} \psi_i^j
\biggl(1-\frac{1}{\rho_i^{2m_i}}\biggr) \left(\frac{1}{2m_i}\rho_i
g_i' + g_i \right)' \frac{\partial}{\partial t_j}=\psi(\pr_{Z({\bold
h})} (\widetilde{\Omega}(\tilde{e}_k\wedge \tilde{e}_{n+1}))).
%\endgathered
$$
The last relations imply that
$$
{\cal L}_{\Bbb R} \{ \widetilde{\Omega} (\tilde{v}_\alpha,
\tilde{v}_\beta)\mid \alpha, \beta \} = V_1+V_2+V_3 \subset {\bold
g}^{4,{\bold h},m, \psi}.
$$
Furthermore,
$$
\pr_{\cal K} ([W_1, V_1] \oplus V_1)={\bold h_0}
\bigoplus\limits_{i=1}^r {\bold h}_i'
$$
by Lemma~2 applied to the manifold~%
$M$.
It remains to observe that
if we choose a~sufficiently general function~%
$h$ then $\pr_{{\cal N}_2}(V_3)={\cal N}_2$, which yields
$$
{\cal L}_{\Bbb R} \{
[\widetilde{\omega}(\tilde{e}_\gamma),\widetilde{\Omega}
(\tilde{v}_\alpha, \tilde{v}_\beta) ], \widetilde{\Omega}
(\tilde{v}_\alpha, \tilde{v}_\beta) \mid  \alpha, \beta, \gamma \} =
{\bold g}^{4,{\bold h},m,\psi}.
$$
Lemma~2 implies that
${\bold{hol}}(N) = {\bold g}^{4,{\bold h},m,\psi}$.
Since the orthogonal coframe
$e^1, \dots, e^n$
is chosen in an~arbitrary neighborhood of~%
$M$, this implies that $\Hol (N)=G^{4,H,m,\psi}$.

\section[]{The Causality Properties of the
Above\--Con\-struc\-ted Metrics}

Recall~[10] that
a~time-oriented Lorentzian manifold~%
$N$
is called {\it globally hyperbolic\/}
whenever it is strongly causal
and for every two points
$p, q, \in N$
the intersection
$J^+(p)
\cap J^-(q)$
of the causal future of~%
$p$
and the causal past of~%
$q$
is a~compact subset of~%
$N$.
The following statement together with Theorem~1
implies the main result of this article.

\vskip0.2cm

{\bf Theorem 2}. {\it The Lorentzian manifolds with the holonomy
groups $G^{1,H}$, $G^{2,H}$, $G^{3,H, \phi}$, and $G^{4,H, m,
\psi}$, constructed in Theorem~{\rm 1}, are globally hyperbolic
for a~suitable choice of~%
$f$,
$A$,
and~%
$\varepsilon$. }

\vskip0.2cm

{\bf Proof}. Consider each type of holonomy groups separately.

\medskip
{\bf Types 1 and 2}.
It is not difficult to observe that
in our construction of the metrics
with the holonomy groups of these types
we can choose as~%
$f$
a~function with compact support in~%
$N$.
Hence,
as
$\varepsilon \rightarrow 0$
the metric~(1) converges to the metric
$
g_0 = 2 d\eta d\xi + ds^2
$
in the fine
$C^0$-%
topology.
It is clear that
 for a~sufficiently small
$\varepsilon>0$
the vector field
$\frac{\partial}{\partial \xi} +
\frac{\partial}{\partial \eta}$
is time-like and orients~%
$N$
in time with respect to both metrics
$g_0$
and~%
$\tilde{g}$.
The space
$(N,g_0)$
is isometric to the direct product of
the~2-dimensional flat Minkowski space
and a~complete Riemannian space~%
$M$;
therefore,
it is a~globally hyperbolic manifold [10, Chapter~2].
Since global hyperbolicity is a~%
$C^0$-%
stable property [10, Chapter~6],
we obtain the required~ result.

\medskip
{\bf Type 3}.
Consider the metric
$
g_0=2 d \eta (d\xi - (1+|\xi|) d\eta) +\frac{1}{2} g.
$
Recall~[10] that
the space of Lorentzian metrics on a~given manifold~%
$N$
is partially ordered as follows:
Assume that
$g_1 \preceq g_2$
($g_1 \prec g_2$)
whenever every light cone of
$g_1$
is
(strictly)
contained in the light cone of
$g_2$
or,
in other words,
for every point
$p \in N$
and every nonzero vector
$X \in T_p N$
the inequality
$g_1(X,X)\leq 0$
implies that
$g_2(X,X) \leq 0$
($g_2(X,X) < 0$).

\vskip0.2cm

{\bf Lemma 3}. {\it $\tilde{g} \preceq g_0$ for a~sufficiently small
$\varepsilon>0$. }

\vskip0.2cm

{\bf Proof}. Take $p \in N$ and $V=V_1+V_2 \in T_p N$, where $V_1$
and $V_2$
are the components of~%
$V$
tangent respectively to
${\Bbb R}^2$
and~%
$M$.
Observe that
$|A(V_2)|^2 \leq C g(V_2,V_2)$.
Then
$$
%\allowdisplaybreaks\gathered
\tilde{g}(V,V) \geq 2d \eta ( d \xi + \varepsilon f d \eta
) (V_1,V_1) - 4 \varepsilon d \eta (V_1) |A (V_2)| +
g(V_2,V_2)
$$
$$
\geq 2d \eta ( d \xi + \varepsilon f d \eta ) (V_1,V_1) - 2
\varepsilon (d \eta (V_1))^2-2 \varepsilon |A (V_2)|^2 + g(V_2,V_2)
$$
$$
\geq 2d \eta ( d \xi + \varepsilon (f-1) d \eta ) (V_1,V_1) +
(1-2 C \varepsilon) g(V_2,V_2)
$$
$$
\geq 2d \eta ( d \xi -(|\xi|+1) d \eta ) (V_1,V_1) +
\frac{1}{2} g(V_2,V_2) \geq g_0(V,V)
%\endgathered
$$
if we choose
$\varepsilon>0$
so small that
$1-2 C \varepsilon
\geq \frac{1}{2}$
and
$$
%\gathered
|\varepsilon (f-1)| = \varepsilon \biggl| h - \xi \sum\limits_{i=1}^r \phi_i
\left(\frac{1}{2m_i}\rho_i g_i' + g_i \right) -1 \biggr|
$$
$$
\leq
|\xi| \biggl|\varepsilon \sum\limits_{i=1}^r \phi_i
\left(\frac{1}{2m_i}\rho_i g_i' + g_i \right) \biggr|+\varepsilon
|h-1| \leq |\xi|+1.
%\endgathered
$$
Consequently,
every non-space-like
(time-like)
vector for~%the metric
$\tilde{g}$
will be non-space-like
(time-like)
for
$g_0$
as well.
The proof of the lemma is complete.

\medskip
It is obvious that
the vector field
$\frac{\partial}{\partial \eta}$
is time-like for
$(N, g_0)$,
and by Lemma~3
for
$(N,
\tilde{g})$
as well.
Hence,
%the field
$\frac{\partial}{\partial \eta}$
defines the direction of time on~%
$N$
with respect to the Lorentzian metrics
$g_0$
and~%
$\tilde{g}$.
Henceforth we will refer to precisely this direction of time.

\vskip0.2cm

{\bf Lemma 4}. {\it The manifold $(N,g_0)$ is stably causal and
globally hyperbolic.}

\vskip0.2cm

{\bf Proof}. To start off,
 consider the metric
$g_1=2 d \eta (d\xi
- (1+|\xi|) d\eta)$
on
${\Bbb R}^2$.
Define the function
$
T=\eta- \ln (|\xi|+2)
$
on ${\Bbb R}^2$ as follows:
If we consider a~non-space-like regular curve
$\gamma (t) = (\xi(t), \eta(t))$
directed into the future
then
$
\dot{\eta} ( \dot{\xi}-(|\xi|+1)\dot{\eta} ) \leq 0.
$
Consequently,
$$
\dot{T}=\dot{\eta}-\frac{\dot{\xi}}{|\xi|+2} \geq
\frac{1}{|\xi|+2} ( (|\xi|+1)\dot{\eta} - \dot{\xi}
)+\frac{\dot{\eta}}{|\xi|+2} > 0.
$$
Hence,~%
$T$
is a~global function of time on
${\Bbb R}^2$,
and so
$({\Bbb R}^2,g_1)$
is stably causal.
It is obvious that~%
$T$
is also a~global function of time on
$(N,g_0)$,
which yields the stable causality of the latter.

Consider now a~pair of points
$p_1=(\xi_1,\eta_1)$
and
$p_2=(\xi_2,\eta_2)$
in
${\Bbb R}^2$.
Direct calculations reveal that
the set
$J^+(p_1) \cap J^-(p_2)$
in
${\Bbb R}^2$
is defined by the inequalities
$$
\eta_1 \leq \eta \leq \eta_2,
\quad
c(\xi_1,\eta_1)+\sgn (\xi) \ln (1+|\xi|) \leq \eta \leq c(\xi_2,\eta_2)+\sgn (\xi) \ln (1+|\xi|),
$$
where
$c(\xi_i,\eta_i) = \eta_i - \sgn (\xi_i) \ln (1+|\xi_i|)$.
It is clear that
 the last inequality defines a~compact set in
${\Bbb R}^2$
independently of the choice of %points
$p_1$
and~%
$p_2$.
Consequently,
$({\Bbb R}^2,g_1)$
is globally hyperbolic.
Hence,
$(N,g_0)$,
which is isometric to the direct product of
$({\Bbb R}^2, g_0)$
and
$(M,g)$,
is globally hyperbolic as well.
The proof of the lemma is complete.

It is not difficult now to obtain the required property.
Lemma 3 implies that
a~function of time for ~
$g_0$
is also a~function of time for
$g$.
Consequently,
$(N,\tilde{g})$
is stably causal,
and so
it is strongly causal.
Furthermore,
by Lemma~4,
if
$p, q \in N$
then
$J^+(p) \cap J^-(q)$
with respect to %the metric
$g_0$
is compact.
Thereby,
the closure of % the set
$J^+(p) \cap J^-(q)$
with respect to~% the metric
$\tilde{g}$
is compact.
It follows from [10, Chapter~3] that
$(N,\tilde{g})$
is globally hyperbolic.

\medskip
{\bf Type 4}.
Take the Riemannian metric
$g'$
on
$M^{n_0} \times C_{m_1} \times \dots \times C_{m_r}$
considered in the proof of Theorem~1.
The Riemannian space
$(M,g)$
is isometric to the product of the flat space
${\Bbb R}^{n-m}$
and the Riemannian manifold
$(M^{n_0} \times C_{m_1} \times
\dots \times C_{m_r},g')$.
Consider
the ``background'' metric on~$N$:
$$
g_0=2 d \eta \left(d\xi - \biggl(1+\biggl|\sum\limits_{k=1}^{n-m} t_k
\biggr|\biggr) d\eta\right) + \sum\limits_{k=1}^{n-m} dt_k^2 +
\frac{1}{2} g'.
$$

\vskip0.2cm

{\bf Lemma 5}. {\it $\tilde{g} \preceq g_0$. }

\vskip0.2cm

{\bf Proof}. Use the notation of the proof of Lemma~3. It is not
difficult to observe that $|A(V_2)|^2\leq C g'(V_2,V_2)$, whence
$$
\tilde{g}(V,V) \geq 2d \eta ( d \xi + \varepsilon (f-1) d
\eta ) (V_1,V_1) + \sum\limits_{k=1}^{n-m} d t_k^2 + (1-2 C
\varepsilon) g'(V_2,V_2) \geq g_0 (V,V)
$$
if we choose
$\varepsilon>0$
so small that
$1-2 C
\varepsilon \geq \frac{1}{2}$
and
$$
%\gathered
|\varepsilon (f-1)| = \varepsilon \biggl| h- \sum\limits_{k=1}^{n-m}
\sum\limits_{i=1}^r \psi_i^k t_k \left(\frac{1}{2m_i}\rho_i g_i' + g_i
\right)-1 \biggr|
$$
$$
\leq
\biggl| \sum\limits_{k=1}^{n-m} t_k \biggr| \biggl|\varepsilon \sum\limits_{i=1}^r
\psi_i^k \left(\frac{1}{2m_i}\rho_i g_i' + g_i \right)
\biggr|+\varepsilon |h-1| \leq \biggl| \sum\limits_{k=1}^{n-m} t_k
\biggr|+1.
%\endgathered
$$
Consequently,
every non-space-like
(time-like)
vector for~%the metric
$\tilde{g}$
will be non-space-like
(time-like)
for
$g_0$
as well.
The proof of the lemma is complete.

As in the previous case,
Lemma~5 enables us to take the field
$\frac{\partial}{\partial \eta}$,
which is time-like for
$(N, g_0)$,
as the time-orienting field for
$(N,\tilde{g})$
and
$(N,g_0)$,
as well as for all auxiliary Lorentzian metrics used below.

To simplify notation,
 put
$$
F=F(\xi, \eta, t_1, \dots,
t_k)=\biggl|\sum\limits_{k=1}^{n-m} t_k \biggr|.
$$
Consider on
${\Bbb R}^{2+n-m}$
the two auxiliary Lorentzian metrics
$$
g_1=2 d \eta (d\xi - (1+F ) d\eta) +
\sum\limits_{k=1}^{n-m} dt_k^2,
$$
$$
g_2=2 \left(d \eta-\delta  d\xi \right) \left(d\xi - \left(1+F +
(1+F)^2 \delta  \right) d\eta\right) + \left( 1-\delta  \right)
\sum\limits_{k=1}^{n-m} dt_k^2,
$$
where
$\delta=\delta(\xi,\eta,t_1,\dots,t_k)$
is a~real function such that
$0 < \delta < 1$
on
${\Bbb R}^{2+n-m}$.

\vskip0.2cm

{\bf Lemma 6}. {\it We have $g_1 \prec g_2$
for a~suitable choice of~%
$\delta$. Moreover, $({\Bbb R}^{2+n-m},g_2)$ is causal, and
consequently $({\Bbb R}^{2+n-m}, g_1)$ is stably causal.}

\vskip0.2cm

{\bf Proof}.
Given a~ vector~%
$V$,
we have
$$
%\gathered
g_1 (V, V)-g_2 (V, V)=
$$
$$
2\delta  \biggl( d\xi^2 +(1+F)^2 d\eta^2 - ( 1 + F + (1+F)^2\delta
) d \xi d\eta +\frac{1}{2} \sum\limits_{k=1}^{n-m} dt_k^2 \biggr)
(V,V) >0
%\endgathered
$$
if we choose~%
$\delta$
so that
$2(1+F(p))\delta(p)<1$
for all
$p \in {\Bbb R}^{2+n-m}$.
Hence,
$g_1(V,V) \leq 0$
implies
$g_2(V,V)<0$
for every nonzero vector~%
$V$,
and thereby
$g_1 \prec g_2$.

Take in $({\Bbb R}^{2+n-m},g_2)$ some closed regular non-space-like
curve \newline $\gamma(s) = (\xi(s), \eta(s), t_1(s),
\dots,t_{n-m}(s))$, for $0\leq s \leq 1$, with
$\gamma(0)=\gamma(1)$, directed into the future. Then
$$
\dot{\eta}-\delta \left( F \right) \dot{\xi} \geq 0
                                \eqno{(5)}
$$
since the time-like vector
$(0,1,0,\dots,0)$
defines the part of the light cone
directed into the future.
Take the closed regular projection
$\gamma_1(s)=(\xi(s),\eta(s))$
of
$\gamma(s)$
onto the plane
${\Bbb R}^2$.
The curve
$\gamma_1 (s)$
can have selfintersections,
but we can always consider a~closed segment of
$\gamma_1(s)$
without selfintersections.
Furthermore,
on this closed segment
the continuity of
$\delta(s) =
\delta(\gamma(s))$
can be violated;
however,
we can deform it into a~continuous function
without changing the extremal values of~%
$\delta$
on
$\gamma(s)$
as follows:
in a~neighborhood of a~discontinuity point
we must make~%the function
$\delta$
vary between the left and right limits at the discontinuity point.
Moreover,~%
(5) will remain fulfilled by linearity.
Therefore,
we assume that
 the closed planar curve
$\gamma_1(s)$
has no selfintersections.

Consider on
${\Bbb R}^2$
the one-dimensional distribution along
$\gamma_1(s)$
given by
$$
\eta-\delta(s) \xi=0.\eqno{(6)}
$$
Since
$\delta<1$,
(6) cannot make a~``full turn'' around
$\gamma_1$.
Consequently,
we can extend~(6) to a~continuous distribution on the whole plane,
which is integrable.
It follows from~(5) that
if
$\gamma_1(s)$
meets transversally one of the integral curves of~(6)
then it will not meet the distribution any more,
which contradicts the closedness of~%
$\gamma_1$.
Thus,
$({\Bbb R}^{2+n-2},g_2)$
contains no closed non-space-like curves,
and consequently
it is causal.
Since
$g_1 \prec g_2$,
we see that
$({\Bbb R}^{2+n-2},g_1)$
is stably causal.
The proof of the lemma is complete.

\vskip0.2cm

{\bf Lemma 7}. {\it The space $({\Bbb R}^{2+n-2},g_1)$ is globally
hyperbolic.}

\vskip0.2cm

{\bf Proof}. By Lemma 6 the space $({\Bbb R}^{2+n-2},g_1)$ is stably
causal, and so it is strongly causal as well. Take $p_i=(\xi_i,
\eta_i, t_{1i}, \dots, t_{(n-m)i}) \in {\Bbb R}^{2+n-m}$ for $i=1,2$
and $p \in J^+(p_1) \cap J^-(p_2)$. Then there exists
a~non-space-like regular curve $ \gamma(s)=(\xi(s),\eta(s), t_1(s),
\dots, t_{n-m}(s)), $ from $p_1$ to $p_2$
through~%
$p$.
Consequently,
the tangent vector~%
$\gamma$
satisfies
$$
\dot{\eta} \geq 0,
\quad
\dot{\xi}-(1+F(s)) \dot{\eta} \leq 0,
\quad
| \dot{F}(s) | \leq (n-m) \sqrt{ 2 \dot{\eta}
(( 1+F(s)) \dot{\eta} -\dot{\xi} ) },
\eqno{(7)}
$$
where
$F(s)=F(\gamma(s))$.
We see that the coordinate~%
$\eta$
is nondecreasing along~%
$\gamma$.
Take the projection
$\gamma_1$
of~%
$\gamma$
onto
${\Bbb R}^2$.
Reparametrize
$\gamma_1$
by the natural parameter~%
$s$
with respect to the standard Euclidean metric on
${\Bbb R}^2$.

Furthermore,
since~%
$\gamma$
is a~smooth curve up to its endpoints,
there exist finitely many closed intervals
subdividing the domain of~%
$\gamma$
so that on every segment either
$(1+F) \dot{\eta} + \dot{\xi} \leq 0$
or
$(1+F) \dot{\eta} +
\dot{\xi} \geq 0$.
Consider firstly an~arbitrary segment
on which
$(1+F) \dot{\eta} + \dot{\xi} \leq 0$.
Then
$\dot{\xi}<0$
and on this segment of the curve
we can consider~
$\xi$
as the parameter,
$\xi_1' \leq \xi \leq \xi_2'$.
It is not difficult to calculate
the maximal slope with respect to the plane
${\Bbb R}^2$
among those generators of the light cone of ~
$g_1$
whose projections onto
${\Bbb R}^2$
satisfy %the inequality
$(1+F)
\dot{\eta} + \dot{\xi} \leq 0$.
Considering that~%
$s$
is the natural parameter for
$\gamma_1$,
we obtain
$$
|\dot{F}(s)|\leq \sqrt{2}(n-m)\frac{\sqrt{2F+2}}{\sqrt{(1+F)^2+1}}.
$$
Hence,
$$
\left| \frac{ dF } {d \xi} \right| \leq \sqrt{2} \left| \frac{d F}
{ds} \right| \leq 2 (n-m)\frac{\sqrt{2F+2}}{\sqrt{(1+F)^2+1}}
$$
on the segment of~%
$\gamma$
under consideration.
Integrating,
we obtain
$
F(\xi) \leq g(\xi)
$
for some function
$g(\xi)\sim |\xi|^{\frac{2}{3}}$
up to multiplying by a~constant
and adding some terms of smaller order of growth as
$|\xi| \rightarrow \infty$.
Now consider on ~
${\Bbb R}^2$
the metric
$
g_3=2 d \eta (d\xi - g(\xi) d\eta).
$
The last inequality implies that
on the segment under consideration
$\gamma_1$
is non-space-like with respect to~%the metric
$g_3$.
However,
by integrating directly the equations for the light rays of %the metric
$g_3$
we verify that
there are two kinds of these rays:
$\eta=\const$
and the curves
with the asymptotics
$\xi \sim \eta^3$.
However,
this means that
$\gamma_1$
cannot leave some bounded domain
$K \subset {\Bbb R}^2$,
which depends only on the initial pair of points
$p_1$
and~%
$p_2$.

Consider one of the segments of the curve
$\gamma$
on which
$(1+F) \dot{\eta}+\dot{\xi} \geq 0$.
Take here
$\eta$
as the parameter.
Rearrange the last of the inequalities in~(6):
$$
\left| \frac{ d F }{d \eta} \right| \leq (n-m) \sqrt{ 2 \left(
\left( 1+F(\xi) \right) -\frac{d \xi}{d\eta} \right) } \leq 2 (n-m)
\sqrt{  1+  F        }.
$$
Integrating,
we obtain
$$
|F(\eta)| \leq g(\eta)
$$
for some function
$g(\eta) \sim \eta^2$
as
$\eta \rightarrow
\infty$.
As in the previous case,
the projection
$\gamma_1$
will be a~non-space-like curve with respect to the metric
$
g_3=2 d \eta \left(d\xi - g(\eta) d\eta\right)
$
on
${\Bbb R}^2$.
Integrating,
we obtain the light rays with the asymptotics
$\eta=\const$
and
$\xi \sim \eta^3$,
which means again the impossibility for %the curve
$\gamma_1$
to leave some bounded domain.

Thus,
the entire curve
$\gamma_1$
cannot leave some bounded domain~%
$K$;
moreover,
the last inequality implies that~%
$F$
is bounded along~%
$\gamma$
by some constant depending only on
$p_1$
and
$p_2$.
Hence,~%
$p$
belongs to some bounded domain in
$R^{2+n-m}$
depending only on
$p_1$
and
$p_2$.
We deduce that
the closure of
$J^+(p_1) \cap J^-(p_2)$
is compact in a~strongly causal space;
consequently,
the space is globally hyperbolic.
The proof of the lemma is complete.

It is not now difficult to finish proving Theorem~2.
The space
$(N,g_0)$
is isometric to the direct product of
the globally hyperbolic space
$({\Bbb R}^{2+n-m},g_1)$
and the complete Riemannian space
$$
\biggl(M^{n_0} \times C_{m_1} \times \dots\times
 C_{m_r},
\frac{1}{2}g'\biggr).
$$
Therefore,
it is globally hyperbolic as well
[10, Chapter~2].
Lemmas~5 and~6 imply that
$(N,\tilde{g})$
is stably causal,
and so
strongly causal.
Moreover,
the closure of the intersection of causal future and causal past
in
$(N,\tilde{g})$
lies within the corresponding intersection in
$(N, g_0)$,
and so it is compact.
It follows from [10, Chapter~3] that
$(N,\tilde{g})$
is globally hyperbolic.

\end{document}